\documentclass[12pt,leqno]{article}
\usepackage[pdftex]{graphicx}
\usepackage{amsfonts}

\newcounter{conjecture}
\newcounter{remark}
\newcounter{corollary}

\newtheorem{theorem}{Theorem}
\newtheorem{lemma}{Lemma}

\newcommand{\lar}{\longrightarrow}

\newcommand{\reals}{\mathbb{R}}
\newcommand{\lll}{\label}

\def \be{\begin{equation}}
\def \ee{\end{equation}}
\def \bt{\begin{theorem}}
\def \et{\end{theorem}}
\def \bea{\begin{eqnarray}}
\def \eea{\end{eqnarray}}
\def \bas{\begin{eqnarray*}}
\def \eas{\end{eqnarray*}}

\def \noi{\noindent}

\def \vski{\vspace{12pt}}
\def \ff{\infty}

\def \({\left(}
\def \){\right)}

\def \intt{\cap}

\def \bc{\begin{center} }
\def \ec{\end{center} }
\def \bs{\begin{slide} }
\def \es{\end{slide} }

\def\square{{\vcenter{\vbox{\hrule height.3pt
         \hbox{\vrule width.3pt height5pt \kern5pt
            \vrule width.3pt}
         \hrule height.3pt}}}}
\def\qed{{\hfill $\square$ \bigskip}}

\begin{document}

\title{Pascal's Hexagon Theorem implies a Butterfly Theorem in the Complex Projective Plane}

\author{Greg Markowsky}

\maketitle

\section{Introduction}

Some time ago I attempted to prove the following for my own entertainment.

\vspace{12pt}

\noindent {\bf Butterfly Theorem.} {\it Let ab be a chord of a circle with midpoint m. Suppose
rs and uv are two other chords that pass through m, as shown below. Let p
and q be the intersections of rv and fs with ab. Then pm = qm.}

\vspace{6pt}
\hspace{.45in}
\includegraphics[width=100mm, height=70mm]{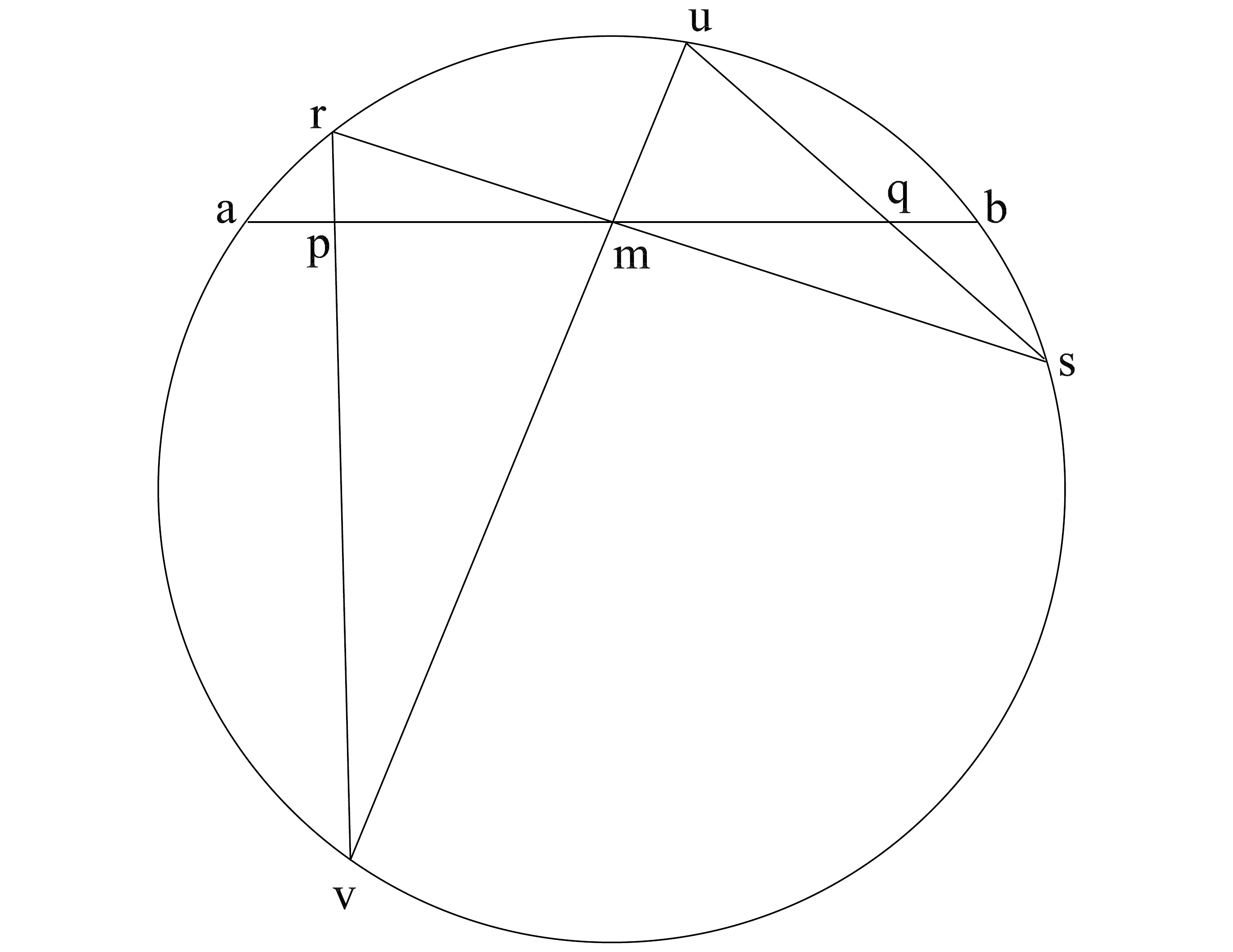}

\noindent Much to my delight, the proof I came up with used Pascal's Hexagon Theorem. A few years later, I learned of the theory of conics in the complex projective plane. I attempted to translate my proof for the circle into this more general setting, and was thrilled to see that it worked very naturally. In what follows, $cr(abcd)$ denotes the cross ratio, and $C$ is a fixed conic in $\mathbb{CP}^2$. Here then is the generalization:

\bt \lll{damn}
Suppose $C$ contains distinct points $a,b$. Let $m$ be any point on $ab$ which is not equal to $a$ or $b$. Let $rs$ and $fg$ be any chords of $C$ which contain $m$. Let $i=rg \intt ab$ and $j=fs \intt ab$. If $p$ is chosen on $ab$ so that $cr(pamb)=-1$, then $cr(pjmi)=-1$ as well.
\et

\includegraphics[width=100mm, height=70mm]{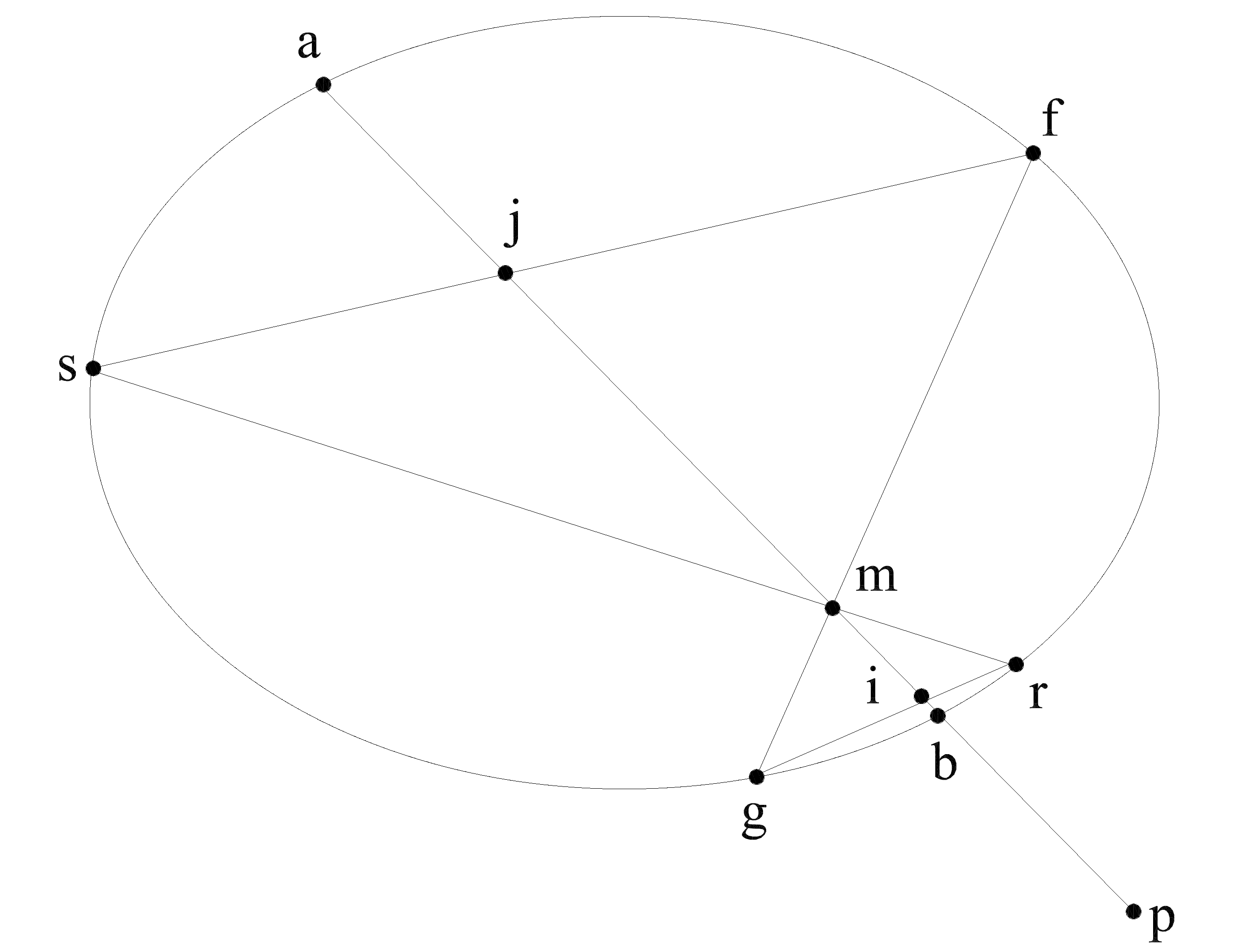}

The next section contains the proof of this theorem, and the final section contains a few notes on the planar case.

\section{Proof of Theorem \ref{damn}}

Our initial aim is to define a natural way to reflect points around any chord of a conic. We need a lemma first.

\begin{lemma} \lll{mono}
Let $k$ be any line which intersects $C$ at two points $u$ and $v$, and let $p$ denote the pole of $k$. Let $l$ be a line through $p$ which intersects $C$ at two points $y,y'$, and let $m$ be a point on $l$. Then $cr(pymy')=-1$ if, and only if, $m=l \intt k$.
\end{lemma}

\includegraphics[width=100mm, height=70mm]{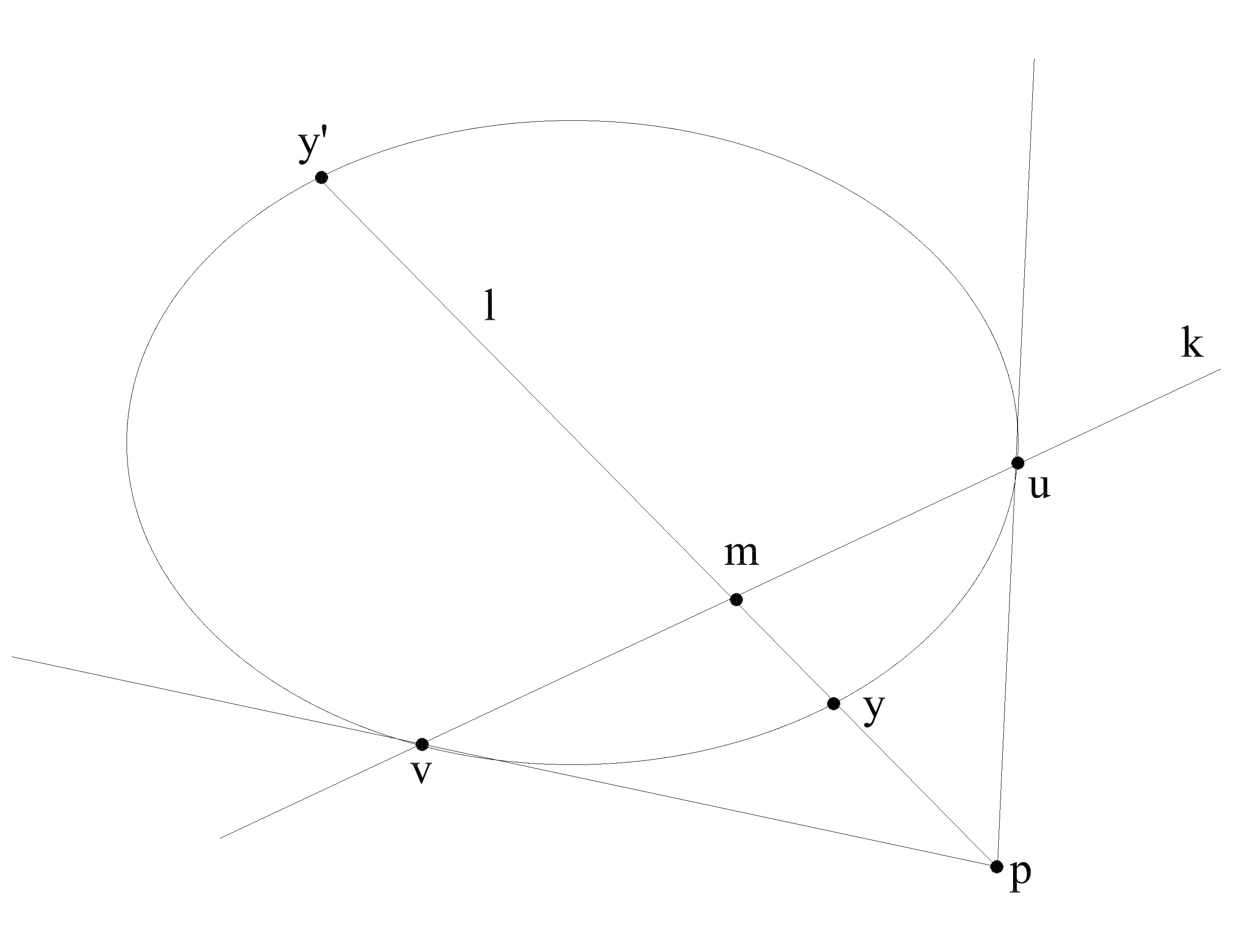}

\noi {\bf Proof:} Choose homogeneous coordinates such that $u=(0:1:0), v=(0:0:1), y'=(1:0:0)$, and $y=(1:1:1)$. Then the equation defining $C$ must be of the form $Axy+Bxz+Cyz=0$ with $A+B+C=0$. We can see that $l:y-z=0$ and $k:x=0$. Let $m=l \intt k$. Then $m=(0:1:1)$. Using the standard formula for the tangent line to a point(see \cite{nambda}), the tangent at $u$ is $Ax+Cz=0$, and the tangent at $v$ is $Cy + Bx=0$. These intersect at the point $p=(1:\frac{-B}{C}:\frac{-A}{C})$. But $p$ also lies on $l$, so that $A=B$. The equation $A+B=-C$ therefore forces $p=(1:1/2:1/2)$. We can projectively map $(p,y,m,y') \lar (1/2,1,\ff,0)$, so we see that indeed $cr(pymy')=-1$. The converse follows easily. \qed

We begin by defining $y$ and $y'$ as referred to in the previous lemma as reflections of each other around $k$. We see that this induces a natural mapping of $C$ to itself. We can extend this to a map for any point $y$ in the plane by $y \lar y'$, where $n=k\intt py$ and $y'$ is the unique point on $py$ such that $cr(pyny')=-1$(this map is defined to be the identity on $k$, and is undefined at $m'$). We will refer to this mapping as the reflection over $k$.

\begin{lemma} \lll{jap}
Let $k$ be a line with pole $p$, let $y$ and $y'$ be points which are reflections of each other over $k$, and let $u$ be a point on $k$. Draw any other line through $p$, and let $t,t'$ be the intersections of this line with $ru,r'u$. Then $t$ and $t'$ are reflections of each other.
\end{lemma}

\includegraphics[width=100mm, height=70mm]{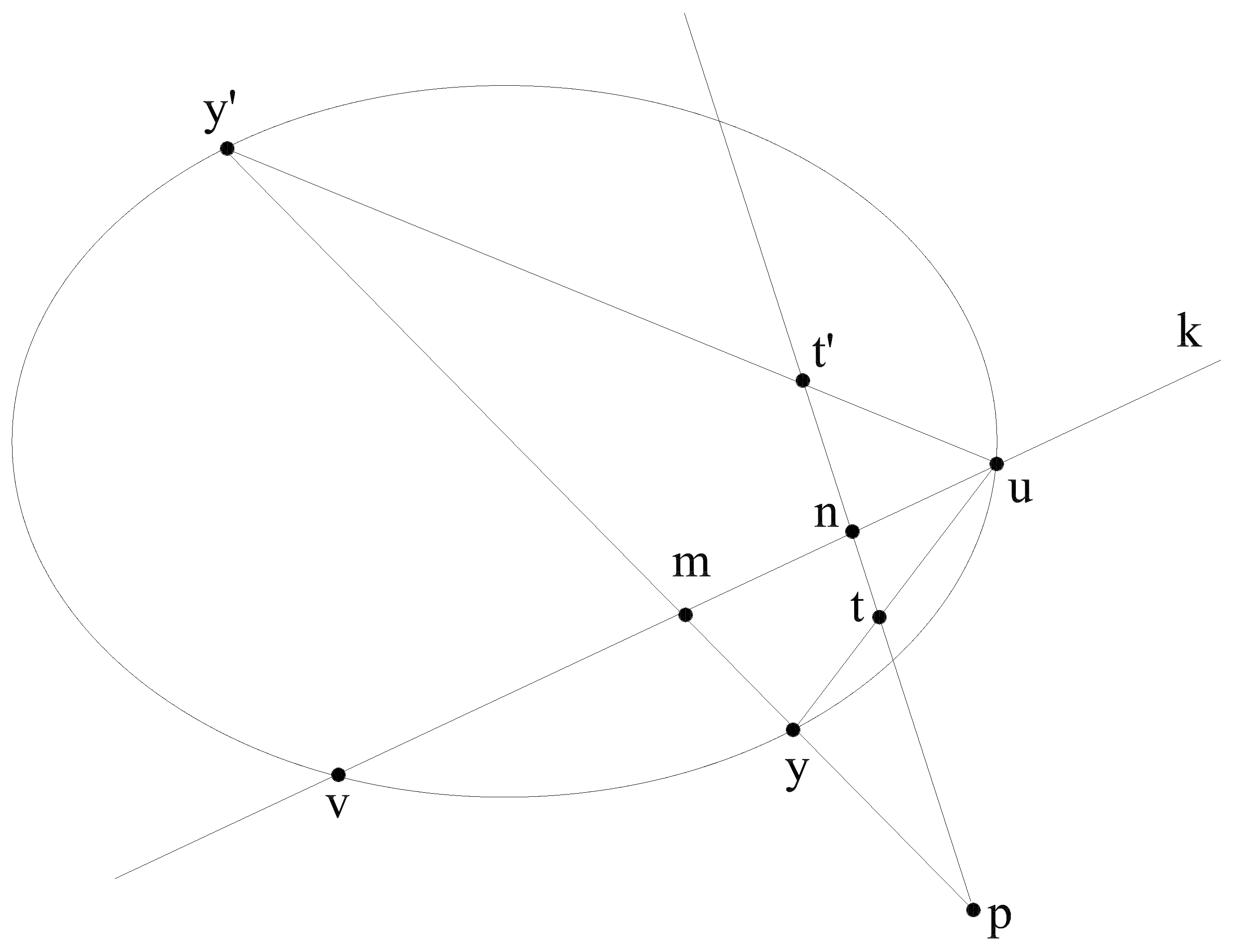}

\noi {\bf Remark:} In the picture it is shown that $u, y, y' \in C$, but this isn't necessary for the lemma to hold.

\vski

\noi {\bf Proof:} Let $m=k \intt py$ and $n=k \intt pt$ . Projecting $pt$ from $u$ onto $py$ we see that $cr(ptnt')=cr(pymy')$, which is $-1$ by Lemma \ref{mono}. \qed

\begin{lemma} \lll{nut}
Suppose $y,y'$ are reflections of each other over a line $k$ with pole $p$, as are $z,z'$. Then $yz \intt y'z' \in k$.
\end{lemma}

\includegraphics[width=100mm, height=70mm]{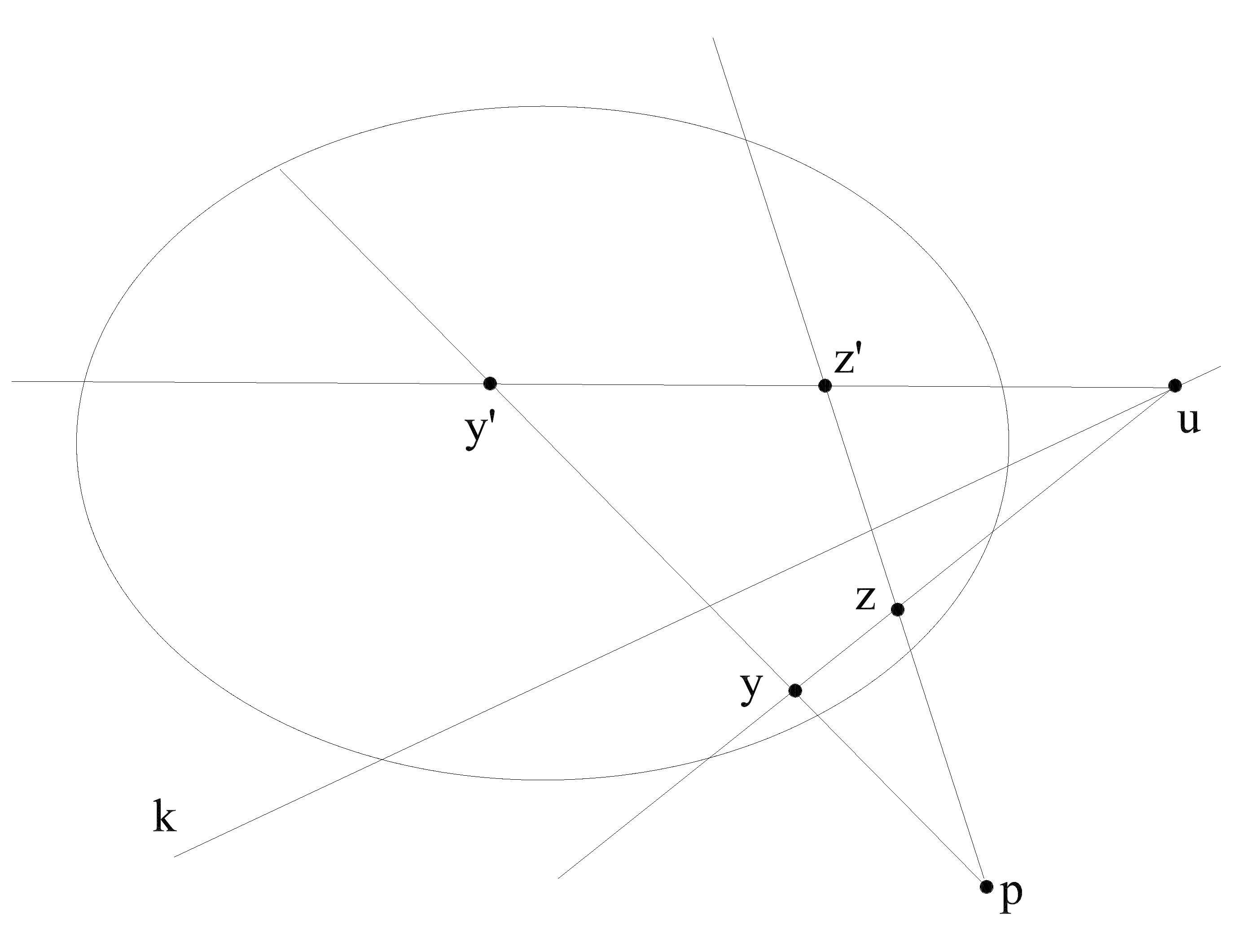}

\noi {\bf Proof:} Let $u=yz \intt k$. Then, by Lemma \ref{jap}, $uy' \intt pz$ is the reflection of $z$, and is therefore
equal to $z'$. Thus, $y'z'$ passes through $u$ as well, and $u=yz \intt y'z' \in k$ \qed

In this situation, we will say that the lines $yz$ and $y'z'$ are reflections of each other.

\begin{lemma} \label{sack}
Let $k$ be a chord containing distinct points $u,v$ on $C$, and let $p$ be the pole of $k$. Let $m$ be a point on $k$, and let $rs$ be a chord of $C$ passing through $m$. Let $r'$ be the reflection of $r$ over $k$. Then $x=r'v \intt su$ lies on $pm$.
\end{lemma}

\includegraphics[width=100mm, height=70mm]{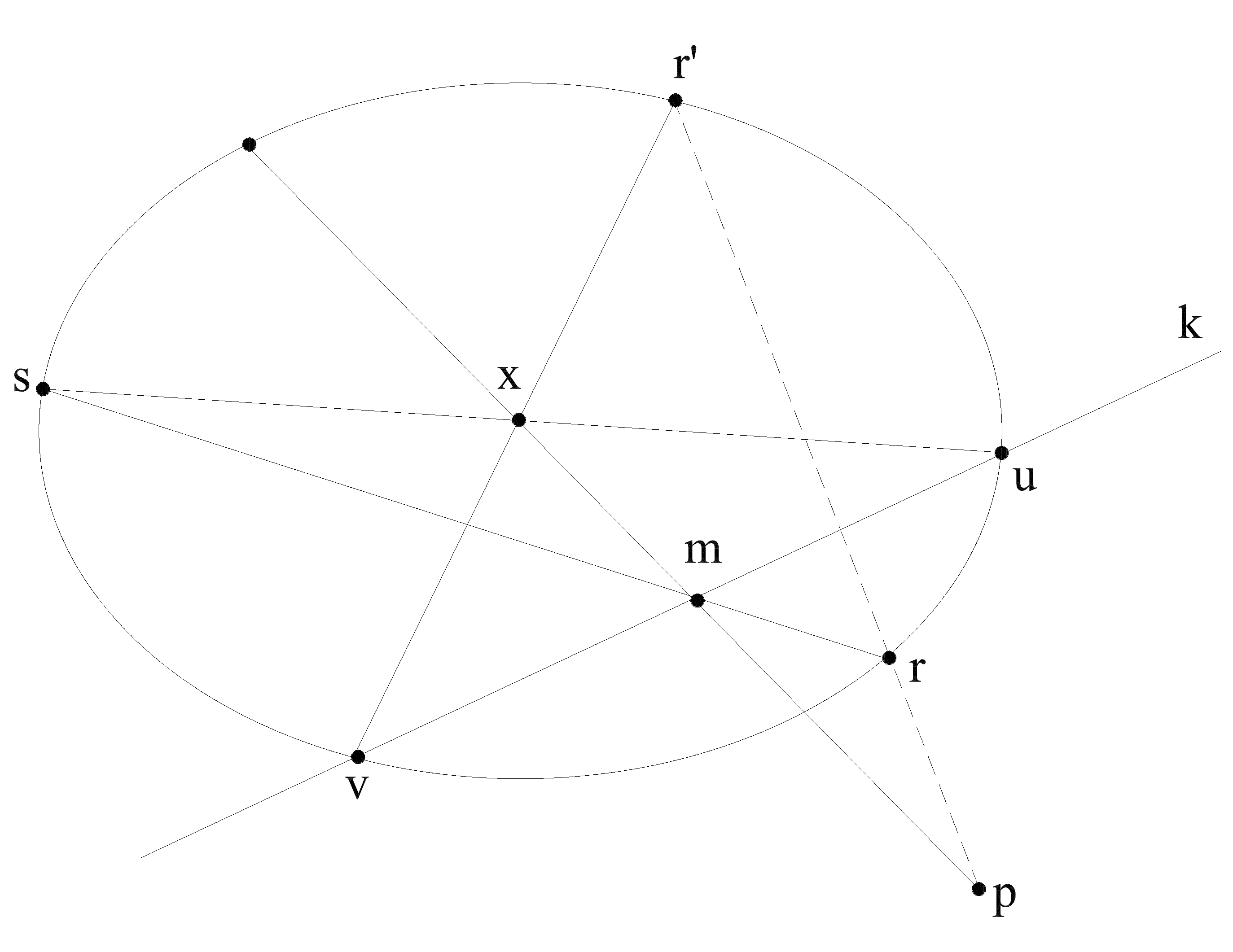}

\noi {\bf Proof:} Let $s'$ be the reflection of $s$, and consider the hexagon $rvr's'us$. Let $z=rv \intt s'u$. By Lemma \ref{sack}, $r's'\intt r's$ lies on $k$, and must therefore be equal to $m$. By Pascal's Theorem, $z,m,$ and $x$ lie on a line. The theorem will be proved if we can show that this line passes through $p$, as is shown in this picture.

\includegraphics[width=100mm, height=70mm]{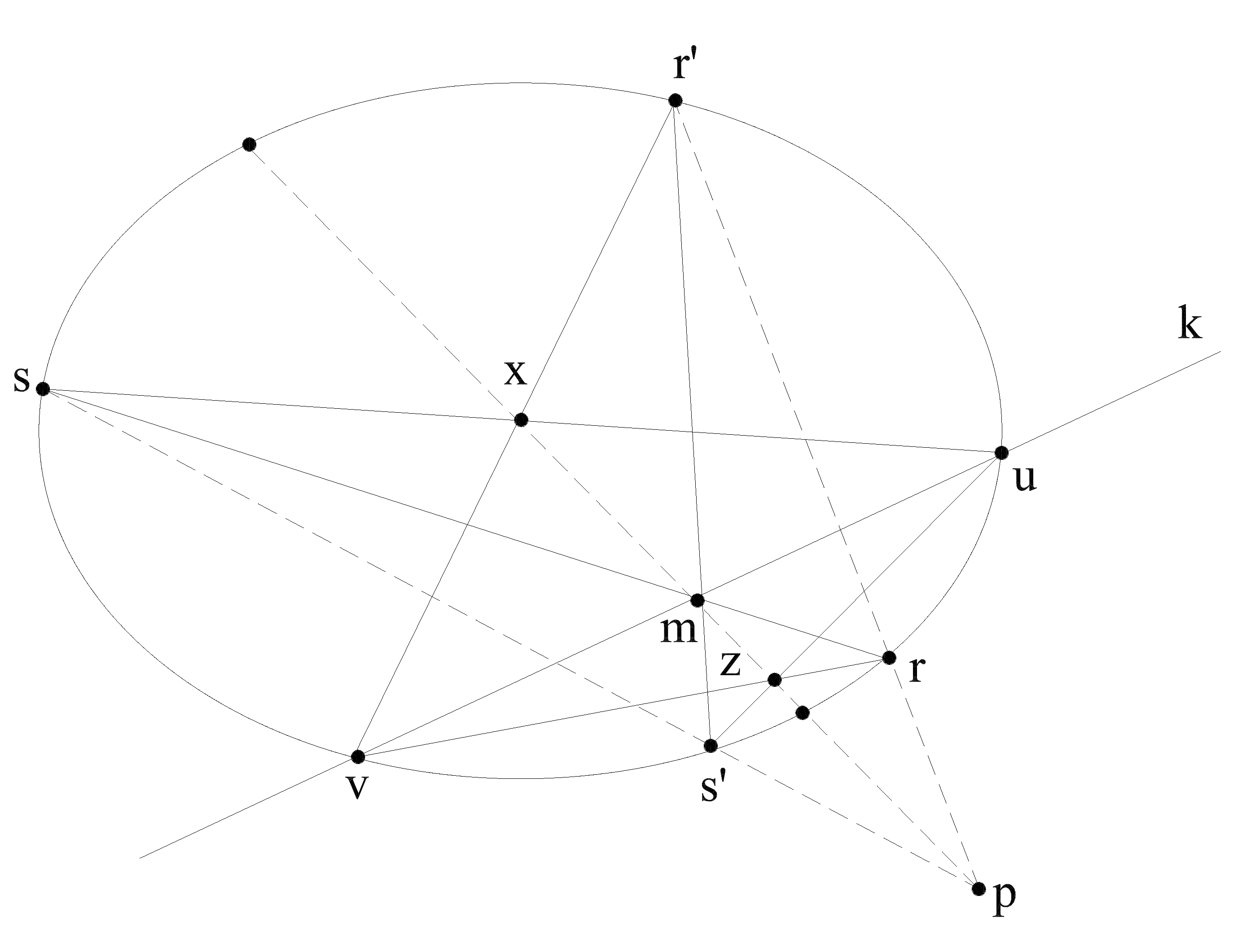}

$su$ and $s'u$ are reflections of each other, as are $rv$, $r'v$. It follows that the reflection of $x$ lies on both $us'$ and $vr$, and is therefore equal to $us' \intt vr$. Thus, $z$ and $x$ are reflections of each other, which implies that $z,x,$ and $p$ lie on a line. Hence, $z,m,x,$ and $p$ are collinear. \qed

\noi {\bf Proof of Theorem \ref{damn}:} In light of what has come before, we need to prove that $i$ and $j$ are reflections of each other. Let $g'$ and $r'$ be the reflections of $g$ and $r$.

\includegraphics[width=100mm, height=70mm]{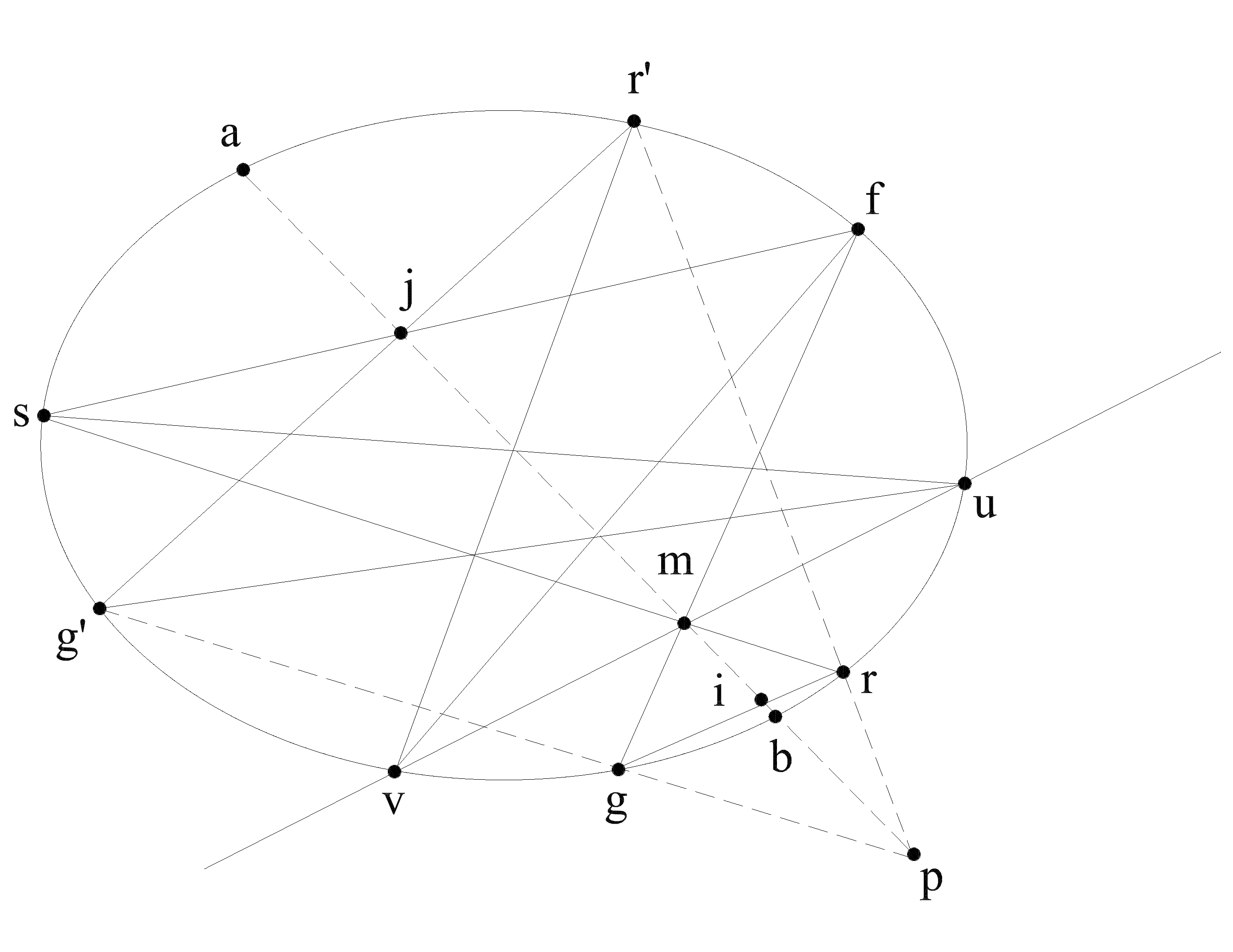}

By Lemma \ref{sack}, $r'v \intt su$ and $fv \intt g'u$ both lie on $pm$. Thus, by Pascal's Theorem applied to the hexagon $r'vfsug'$, $r'g' \intt sf=j$ lies on $pm$ as well. But $r'g'$ is the reflection of $rg$, so by Lemma \ref{nut}, $rg \intt pm=i$ is the reflection of $r'g' \intt pm=j$, and we are done. \qed

\section{Remarks on the planar case}

In \cite{markow}, this method of proof is used to deduce the Butterfly Theorem for conics in $\reals ^2$, with one exception. The case in which the initial chord intersects a hyperbola once on each of the branches of the hyperbola could not be covered while staying entirely in $\reals^2$, since the relevant polar in that case did not intersect the hyperbola at any point. To get around this difficulty, we consider $\mathbb{RP}^2$ as embedded in $\mathbb{CP}^2$ as the set of fixed points of the map $(z_1:z_2:z_3)\lar(\bar{z_1}:\bar{z_2}:\bar{z_3})$. In this larger space, all lines intersect the conic, and we arrive at no difficulties. Therefore, the following generalization of the Butterfly Theorem in $\reals ^2$ is obtained as a corollary to the above work:

\bt \lll{cutl}
Let $C$ be a conic in the plane. Let a point $m$ be on a chord intersecting $C$ at two distinct points $a$ and $b$. Let $rs$ and $uv$ be two chords passing through $m$. Let $p$ and $q$ be the intersections of $ru$ and $sv$ with $ab$. Let $m'$ be the unique point(possibly $\ff$) on $ab$ so that $cr(m'amb)=-1$. Then $cr(m'pmq)=-1$ as well.
\et

\noi {\bf Proof:} Suppose $C$ is given by $Ax^2 + By^2 +Qxy +Dx +Ey +F=0$. Then $Ax^2 + By^2 +Qxy +Dxz +Eyz +Fz^2=0$ gives the extension of $C$ to $\mathbb{CP}^2$. Since $m,a,b \in \mathbb{RP}^2$, $m' \in \mathbb{RP}^2$ as well. It follows as above that $p$ and $q$ are reflections of each other over the polar of $m'$. This polar also passes through $m$, though it does not necessarily intersect $C$ in $\mathbb{RP}^2$. Whether or not the polar intersects $C$ in $\mathbb{RP}^2$, we have $cr(m'pmq)=-1$, and we are done. \qed

\section{Acknowledgement}

I would like to thank Jihun Park for several helpful conversations.

\def\noopsort#1{} \def\printfirst#1#2{#1} \def\singleletter#1{#1}
   \def\switchargs#1#2{#2#1} \def\bibsameauth{\leavevmode\vrule height .1ex
   depth 0pt width 2.3em\relax\,}
\makeatletter \renewcommand{\@biblabel}[1]{\hfill#1.}\makeatother

\it{Greg Markowsky}

\it{gmarkowsky@gmail.com}

\end{document}